\documentclass[aps,pra,preprint,amsmath,amsfonts,amssymb,a4paper,eqsecnum,superscriptaddress]{revtex4-1}
\usepackage {eucal}
\newtheorem{definition}{Definition}[section] 

\newtheorem{remark}{Remark}

\newcommand{\pai}[2]{\langle\,#1\,,#2\,\rangle}  
\newcommand{\map}[3]{#1\colon#2\rightarrow#3}  
\renewcommand{\ddot}[1]{\frac{d^2 #1}{dt^2}}  
\renewcommand{\dddot}[1]{\frac{d^3 #1}{dt^3}}  
\renewcommand{\ddddot}[1]{\frac{d^4 #1}{dt^4}}

\renewcommand{\sin}{\mathrm{sin}}
\renewcommand{\cos}{\mathrm{cos}}
\renewcommand{\tan}{\mathrm{tan}}
\newcommand{\id}{\mathrm{id}}

\newcommand{\at}[1]{\Big\vert_{#1}}

\newcommand{\set}[2]{\left\{\,#1\left.\vphantom{#1#2}\,\right\vert\,#2\,\right\}} 

\def\T{\CMcal T}

\let\X\baseX
\let\V\baseV
\let\P\baseP

\newcommand{\bra}[3][]{\lbrack\! \lbrack#2,#3\rbrack\! \rbrack^{#1}}

\newcommand{\qquand}{\qquad\mathrm{\myand}\qquad}

\def\myand{and}

\def\D{\mathcal{D}}
\newcommand{\TED}[1][]{\CMcal{T}^{E}_{#1}\D}
\newcommand{\TDD}[1][]{\CMcal{T}^{\D}_{#1}\D}

\newcommand{\cf}{\boldsymbol{\kappa}} 

\newcommand{\cinfty}[1]{C^\infty(#1)}
\def\pd#1#2{\displaystyle{\frac{\partial#1}{\partial#2}}}
\usepackage{xy, color}
\xyoption{all}

\def\Real{\mathbb{R}}

\begin{document}

\title{Some applications of quasi-velocities in optimal control}

\author{L\'\i gia Abrunheiro}
\affiliation{
ISCA - Universidade de Aveiro, 
Rua Associa\c c\~ao Humanit\'aria dos Bombeiros de Aveiro, 
3811-902 Aveiro - Portugal 
\email{abrunheiroligia@ua.pt}
}
\author{Margarida Camarinha}
\affiliation{Departamento de Matematica, 
Universidade de Coimbra,
Largo D. Dinis, Apartado 3008,
3001-454 Coimbra
Portugal 
\email{mmlsc@mat.uc.pt}
}
\author{Jos\'e F.\ Cari\~nena}
\affiliation{Departamento de F\'\i sica Te\'orica 
Universidad de Zaragoza,
50009 Zaragoza, Spain
\email{jfc@unizar.es}
}
\author{Jes\'us Clemente-Gallardo}
 \affiliation{ BIFI-Departamento de F\'\i sica Te\'orica 
Edificio I+D, Campus R\'{\i}o Ebro, Universidad de Zaragoza, 
C/Mariano Esquillor s/n, 50018 Zaragoza, Spain
\email{jesus.clementegallardo@bifi.es}
}
\author{Eduardo Mart\'{\i}nez}
\affiliation{Departamento de Matem\'atica Aplicada, 
Facultad de Ciencias,
Universidad de Zaragoza,
50009 Zaragoza, Spain
\email{emf@unizar.es}
}
\author{Patricia Santos}
\affiliation{CMUC- Universidade de Coimbra and 
Instituto Superior de Engenharia de Coimbra,
Rua Pedro Nunes - Quinta da Nora
3030-199 Coimbra, Portugal
\email{patricia@isec.pt}}



\begin{abstract}
  In this paper we study optimal control problems for nonholonomic
  systems defined on Lie algebroids by using quasi-velocities. We
  consider both kinematic, i.e. systems whose cost functional depends
  only on position and velocities, and dynamic optimal control
  problems, i.e. systems whose cost functional depends also on
  accelerations. The formulation of the problem directly at the level
  of Lie algebroids turns out to be the correct framework to explain
  in detail similar results appeared recently
  \cite{maruskin:2007p3969}.
  We also provide several examples to illustrate our construction.
\end{abstract}

\maketitle
\section{Introduction}
The principles of analytical mechanics established by D'Alembert, Lagrange,
Gauss and Hamilton can also be contemplated from additional mathematical
perspectives providing us methods for understanding Nature's law from new
viewpoints
which may be helpful in solving specific problems and clarifying the way in
which Nature behaves. The traditional techniques were only applied to very simple
models but current technology needs efficient algorithms in areas ranging from
robotics to spacecraft design. Furthermore the computer development with the
corresponding capability of computation suggests the convenience of analysing
different formulations to yield the differential equations  for
multibody dynamics that involve a certain number of constraints.

There exist different techniques to deal with such constrained systems. The
geometric framework of manifolds replacing Euclidean spaces allows us to give
a formulation  for systems with holonomic constraints in terms of generalised
coordinates and free of  Lagrange multipliers. However it is not clear how to
choose generalised coordinates improving computational efficiency: Kane's
method \cite{Kane:1983p3963,Kane:1985p3955} or the Maggi equations formulation 
\cite{Chen:2008p3084}. 
Another recent alternative formulation is given in \cite{Kalaba:2001p4196}.

Nonholonomic constraints are very relevant and appear in many problems in
physics and engineering, and in particular in control theory.  Such
nonholonomic  constraints restrict possible virtual displacements and
when taking into account such constraints d'Alembert-Lagrange
principle leads to Boltzmann--Hamel equations
\cite{Cortes:2006p2369,Cortes:2004p3190,maruskin:2007p3969,MaruskinBloch:2010p2554,
Papastavridis:1995p3990}. 

The concept of quasi-velocity (or generalised velocity) \cite{Carinena:2007p2820}
is of a great relevance in the study of mechanical systems, 
mainly in nonholonomic ones because  the conditions of nonholonomic
constraints can be expressed in a  simpler form.
 Boltzmann-Hamel equations, Gibbs--Appell and Gauss principles, for
 instance,   make use of quasi-coordinates (also called nonholonomic
 coordinates) and the Hamel symbols \cite{Koiller:1992p4308}. The use of
 quasi-velocities in the dynamics of 
nonholonomic system with symmetry has recently been investigated under 
different approaches in
\cite{Bloch:2009p2783,Cortes:2006p2369,Cortes:2009p2371,Cortes:2004p3190} 
and Hamel's equations have been recovered from this perspective.

It has been shown in recent papers \cite{Carinena:2007p2820,
  Cortes:2006p2369,Cortes:2009p2371} that 
the appropriate geometric framework for studying systems with linear
nonholonomic constraints is the framework of Lie algebroids. The
geometric approach to mechanics uses tangent bundles in the Lagrangian
formulation and tangent bundles are
but  particular instances of algebroids. The usual geometric approach to
Lagrangian formalism was then developed in this extended framework of Lie
algebroids \cite{Carinena:2001p2844,Martinez:2001p1130,Martinez:2001p4348,
  Weinstein:1996p4160}  the main advantage
being that such structure arises in 
reduction processes from tangent bundles when the vertical endomorphism
character  is not projectable \cite{{Carinena:2007p2822}}. The geometrical
construction (see \cite{Martinez:2001p1130}) is based on the generalization of the
usual symplectic description of Lagrangian (or Hamiltonian) mechanics
on tangent (or cotangent) bundles.  The dynamics is then defined by a
function of a Lie algebroid (for the Lagrangian formalism) or its dual
(for the Hamiltonian one). Considering, for the sake of simplicity,
only the Lagrangian case now,  the solutions of
the dynamics correspond to 
the analogue of the concept of SODE (considered as the section of the
 bundle $T(TM)$ for a Lagrangian defined on $TM$). But  in the
Lie algebroid case the concept of second tangent bundle is subtle to
define and the notion of prolongation of a Lie algebroid is necessary:

\begin{definition}
Let $(E, \bra{\ }{\ }, \rho)$ be a Lie
algebroid ($\tau:E\to M$) over a manifold $M$ and $\nu: P \to M$ be a
fibration. For every point $p\in P$ we consider the vector space
$$
\T_{p}^EP =\set{(b,v)\in E_x\times T_pP}{\rho(b)=T_p\nu(v)},
$$
where $T\nu: TP \to TM$ is the tangent map to $\nu$ and
$\nu(p)=x$ . The set $\T^{E}P=\bigcup_{p\in P} \T_{p}^E P$ has a natural vector
bundle structure over $P$, the vector bundle projection $\tau^E_P$ being just the
projection $\tau^E_P(b,v)=\tau_P(v)$.

The vector bundle $\map{\tau^E_P}{\T^E P}{P}$ can be endowed with a Lie
algebroid structure (see \cite{Martinez:2001p1130}). The Lie algebroid $\T^E P$ is
called the \textbf{prolongation of }$\map{\nu}{P}{M}$ with respect to
$E$ or the ${E}$-tangent bundle to $\nu$. 
\end{definition}

When we consider the case
$P=E$, the resulting entension generalizes the notion of second
tangent bundle to the Lie algebroid framework and the sections of this
bundle represent the analogue of SODEs in the usual case.

By using this notion,  the symplectic formalism of Lagrangian systems 
in geometric mechanics is easily extended to this more
general setting. In a similar way, Hamiltonian  mechanics and
discrete mechanics can be also generalized to the new
framework.

It is also well known from geometric control theory that many of the
techniques of classical mechanics can be used in control theory 
\cite{Cortes:2006p2369,Cortes:2004p3190,Martinez:2004p4316}. The
theory of Lie algebroids can be 
applied to deal with control problems and  the application of such geometric
tools is very useful for a better understanding of different control
problems. This is our motivation for developing the theory of optimal control
theory using the properties of Lie algebroid theory which is going to be the
appropriate approach to Boltzmann--Hamel equations.

Throughout the paper, we consider a Lie algebroid
$\tau:E\to M$ with anchor mapping $\rho:E\to TM$. When coordinates are
required, we consider a local basis $\{x^i\}$ for the base
manifold $M$, and a basis of sections $\{e_{\alpha}\}$ for the bundle
$E$ which provides a set of coordinates $(x^i, y^\alpha)$ for the Lie
algebroid. The anchor mapping is then represented by the set of
functions $\{ \rho_i^\alpha\}$ and the Lie algebra structure by the
structure functions $\{ C_{\alpha\beta}^\gamma \}$. When considering
the dual bundle $E^*$, the dual basis of 
sections  $\{e^{\alpha}\}$  is chosen and the corresponding
coordinates are denoted as  $(x^i, \mu_\alpha)$.  The corresponding
exterior differential $d:\mathrm{Sec}\bigwedge^kE^*\to
\mathrm{Sec}\bigwedge^{k+1}E^*$  defines the corresponding
algebroid cohomology, with respect to which the concept of symplectic
or presymplectic form can be defined. This concept will be used later
on when providing the geometrical framework for the maximum principle.

Finally, we consider the coordinate functions for the extensions
of the Lie algebroid $E$ by a bundle $P$. We consider the general case
although in practice we use only the case $P=E$ and the
case $P=\D\subset E$, for $\D$ a subbundle of the Lie algebroid. In
any case, considering local coordinates
$(x^i,u^\beta)$ on $P$ and a local basis $\{e_\alpha\}$ of sections of
$E$, we can define a local basis $\{\X_\alpha,\V_\beta\}$ of sections of 
$\T^E P$ by
$$
\X_\alpha(p)
=\Bigl(p,e_\alpha(\nu(p)),\rho^i_\alpha\pd{}{x^i}\at{p}\Bigr) \qquand
\V_\beta(p) = \Bigl(p,0,\pd{}{u^\beta}\at{p}\Bigr).
$$
If $z=(p,b,v)$ is an element of $\T^E P$, with $b=z^\alpha
e_\alpha$, then $v$ is of the form $v=\rho^i_\alpha
z^\alpha\pd{}{x^i}+v^\beta\pd{}{u^\beta}$,  
and  we can write   
$z=z^\alpha\X_\alpha(p)+v^\beta\V_\beta(p)
$.Vertical sections are linear combinations of $\{\V_A\}$.
Analogously, when considering the dual object, $\T^E P^*$, we use a local
basis $\{ \X_\alpha,\P_\beta\}$, where $\{\P_\beta\}$ are the sections
corresponding to the vertical elements  $\P_\beta(p) =
\Bigl(p,0,\pd{}{\nu_\beta}\at{p}\Bigr)$, for $\{(x^i,\nu_\beta)\}$
being a set of coordinates for the bundle  $P^*$. 

This article is organized in the following way. 
We address the interested reader to
\cite{Cortes:2006p2369,Cortes:2004p3190,Martinez:2004p4316} for a detailed description of
Lie algebroids and the construction of general optimal control
systems. We provide a short review of those results in 
 Sections \ref{control} to \ref{sec:dynam-optim-contr} while
 incorporating the formalism of quasi-velocities on them:  Section
 \ref{control} deals with  the optimal  control theory and the
 Pontryagin maximum principle \cite{Martinez:2004p4316},  kinematic optimal control
 is studied in Section  \ref{sec:kinem-optim-contr} and dynamical
 aspects are the aim on Section \ref{sec:dynam-optim-contr}. The
 theory is illustrated in Section \ref{sec:examples} with several
 examples.

\section{Optimal control theory}
\label{control}

As it is well known, optimal control theory is a generalization of classical
mechanics.  The central result in the theory of optimal
control systems is Pontryagin maximum principle. The reduction of optimal
control problems can be performed within the framework of Lie algebroids,
see~\cite{Martinez:2004p4316}. This was done as in the case of classical mechanics, by
introducing a general principle for any Lie algebroid and later studying the
behavior under morphisms of Lie algebroids. See also~\cite{Grabowski:2009p3428} for a
recent direct proof of Pontryagin principle in the context of general
algebroids.

\subsubsection*{Pontryagin maximum principle
~\cite{Martinez:2004p4316}} 

By a control system on a Lie algebroid $\map{\tau}{E}{M}$ with control space
$\map{\pi}{B}{M}$ we mean a section of $\sigma:B\to E$ along $\pi$. A
trajectory of the system 
$\sigma$ is an integral curve of the vector field $\rho \circ \sigma$ along $\pi$.
\[
\xymatrix{%
&\ar[d]^\tau E\ar[r]^\rho& TM\ar[ld]^{\tau_M}\\
B\ar[ru]^{\sigma}\ar[r]_\pi&M&&
}
\]

Given a cost function $L\in\cinfty{B}$ we want to minimize the integral of $L$
over some set of trajectories of the system satisfying some boundary
conditions. Then we define the Hamiltonian function $H\in\cinfty{E^*\times_MB}$ by
$H(\mu,u)=\langle\mu,\sigma(u)\rangle -L(u)$ and the associated Hamiltonian control system  
$\sigma_H$ (a section along $\map{\mathrm{pr}_1}{E^*\times_MB}{E^*}$ of $\T^E{E^*}$ ) defined 
on a subset of the manifold $E^*\times_M B$, by means of the symplectic equation
\begin{equation}
\label{pontryagin}
i_{\sigma_H}\Omega=dH,
\end{equation}
where $\Omega$ is the canonical symplectic form defined on the bundle
$E$ (i.e. a section of the bundle $\bigwedge ^2E^*$ which is closed
for the differential calculus defined on the Lie algebroid).
The integral curves of the vector field $\rho(\sigma_H)$ are said to be the critical
trajectories.

In the above expression, the meaning of $i_{\sigma_H}$ is as follows:
Let 
$\map{\Phi}{E}{E'}$ be a  morphism of the bundle $T\to M$ over a map
$\map{\varphi}{M}{M'}$ and let $\eta$ be a 
section of $E'$ along $\varphi$. If $\omega$ a section of $\bigwedge^pE'{}^*$ then $i_\eta\omega$ 
is the section of $\bigwedge^{p-1}E^*$ given by
$
(i_\eta\omega)_m(a_1,\ldots,a_{p-1})=\omega_{\varphi(n)}(\eta(m),\Phi(a_1),\ldots,\Phi(a_{p-1}))
$
for every $m\in M$ and $a_1,\ldots,a_{p-1}\in E_m$. In our case, the map $\Phi$ is the
prolongation $\map{\T{\mathrm{pr}_1}}{\T^E{(E^*\times_MB)}}{\T^E{E^*}}$ of the
map $\map{\mathrm{pr}_1}{E^*\times_MB}{E^*}$ (this last map fibered over the identity in
$M$), and $\sigma_H$ is a section along $\mathrm{pr}_1$. Therefore, 
$i_{\sigma_H}\Omega-dH$ is a section of the dual bundle to $\T^E{(E^*\times_MB)}$.

It is easy to see that the symplectic equation~(\ref{pontryagin}) has a unique
solution defined on the following subset
\[
S_H=\set{(\mu,u)\in E^*\times_M B}{\pai{dH(\mu,u)}{V}=0\mathrm{\,\,
    for\,\,\, all \,\,\,}V\in\mathrm{Ker}\,\T{\mathrm{pr}_1}}. 
\]
Therefore, it is necessary to perform a stabilization constraint algorithm to
determine the submanifold where  integral curves of $\sigma_H$ do exists.

In local coordinates, the solution to the above symplectic equation is
\[
\sigma_H=\pd{H}{\mu_\alpha}\X_\alpha-\left[\rho^i_\alpha\pd{H}{x^i}+
\mu_\gamma C^\gamma_{\alpha\beta}\pd{H}{\mu_\beta}\right]\P^\alpha,
\]
defined on the subset where $\displaystyle{ \pd{H}{u^A}=0}$,
and therefore the critical trajectories are the solution of the
differential-algebraic equations
\begin{equation}
\label{pontryagin.coordinates}
\dot{x}^i=\rho^i_\alpha\pd{H}{\mu_\alpha}; \quad
\dot{\mu}_\alpha
=-\left[\rho^i_\alpha\pd{H}{x^i}+\mu_\gamma C^\gamma_{\alpha\beta}\pd{H}{\mu_\beta} \right];
\quad
0=\pd{H}{u^A}.
\end{equation}
Notice that $\displaystyle{\pd{H}{\mu_\alpha}}=\sigma^\alpha$.

One can easily see that whenever it is possible to write $\mu_\alpha=p_i\rho^i_\alpha$ then
the above differential equations reduce to the critical equations for the
control system $Y=\rho(\sigma)$ on $TM$ and the function $L$. Nevertheless it is not
warranted that $\mu$ is of that form. For instance in the case of a Lie algebra,
the anchor vanishes, $\rho=0$, so that the factorization $\mu_\alpha=p_i\rho^i_\alpha$ will not
be possible in general.

\section{Kinematic Optimal Control}
\label{sec:kinem-optim-contr}
Let  $\tau:E\to M$ be a Lie algebroid and $\D$  a constraint
distribution. Given a cost function $\map{\cf}{E}{\Real}$, we consider
the following kinematic optimal control problem: we can 
control directly all the (constrained) velocities, and we want to minimize some
cost functional
\[
I(a)=\int_\alpha^\beta \cf(a(t))\,dt,
\]
for $a:[\alpha,\beta]\subset\Real\to E$ over the set of admissible curves taking values
in $\D$. We use coordinates $(x^i, y^a)$ to denote the elements of
this bundle, where $y^a$ will represent the coordinates with respect to some
basis of section for $\D$, as in the last section. We use
capital indices $A, B, C, \ldots$ to represent the coordinates $\{y^A,
y^B, \ldots \} $ corresponding to the elements of the fiber of $E$ not
contained in $\mathcal{D}$. Analogously, $\{\mu^A, \mu^B,\ldots \}$
represent the coordinates for the fiber elements in $E^*$
not corresponding to sections dual to the elements in $\mathcal{D}$.

\begin{remark}

Whenever the cost function $\cf$ is a quadratic function defined on $\D$, the
problem that we are considering is just the problem of sub-Riemannian
geometry. In the case of a degree~1 homogeneous cost function this is
sub-Finslerian geometry, and in the more general case this  problem can be
called sub-Lagrangian problem.

\end{remark}

Since we can control directly the velocities or pseudovelocities, the
control bundle is $B=\D$ and the system map $\sigma:\D\to E$ is just
the canonical inclusion $\sigma(a)=a$. 
\[
\xymatrix{
& E \ar[d]_\tau \ar[r]^{\rho}& TM \ar[dl]^{\tau_M}\\
            \D\ar[ur]^{\sigma }\ar[r]& M&
}
\]
Pontryagin Hamiltonian is a function
$H\in C^\infty(E^*\times_M B)$ defined as
$H(\mu,b)=\langle\mu,\sigma(b)\rangle-\cf(b)$, which in coordinates
reads 

\begin{equation}
  \label{eq:hamiltonian}
H(x^i,\mu_a,\mu_A,u^a)=\mu_au^a-\cf(x^i,u^a).
\end{equation}

The Maximum principle imposes the choice of the control functions such that
\begin{equation}
  \label{eq:optimal}
  \mu_a=\frac{\partial \cf}{\partial u^a}\qquad\left(\mathrm{from
      \,\,\frac{\partial H}{\partial u^a}=0}\right). 
\end{equation}
Under appropriate regularity conditions the set $S_H$ of solutions of this
equation is a submanifold of $E^*\times_M B$, which we call the critical
submanifold. Frequently, this set is but the image of a section of $E^*\times_M B\to
E^*$,  locally given by  
\begin{equation}
  \label{eq:ua}
u^a=u^a(x, \mu).
\end{equation}
Thus the set $S_H$ is diffeomorphic to $E^*$ and the optimal Hamiltonian, with
local expression $H(x^i,\mu_\alpha,u^a(x,\mu))$, defines via the canonical symplectic
form a Hamiltonian system on $E^*$. 

The restriction of the Pontryagin-Hamilton equations to this submanifold
provides us with the control system
\begin{equation}
  \label{eq:controlsys}
   \dot{x}^i=\rho^i_au^a
\end{equation}
and the dynamics of the momenta
\begin{eqnarray}
\label{eq:dotm-=-leftrh}
\dot{\mu}_a
&=-\left[\rho^i_a\pd{H}{x^i}+\mu_c C^c_{ab}u^b+\mu_B C^B_{ab}u^b
\right]\\
\dot{\mu}_A
&=-\left[\rho^i_A\pd{H}{x^i}+\mu_c C^c_{Ab}u^b+\mu_B C^B_{Ab}u^b
\right].
\end{eqnarray}
These are the equations to be solved and to be used to determine, by using the mapping
(\ref{eq:ua}), the control functions defining the solution optimizing the
value of the cost function.
By substitution of $\mu_a=\displaystyle{\pd{\cf}{u^a}}$ into these equations, and taking into
account that $\pd{H}{x^i}=-\displaystyle{\pd{\cf}{x^i}}$, we get
\begin{eqnarray}
\dot{x}^i&=\rho^i_au^a; \qquad
\frac{d}{dt}\left(\pd{\cf}{u^a}\right)-\rho^i_a\pd{\cf}{x^i}+\pd{\cf}{u^c}
C^c_{ab}u^b+\mu_B C^B_{ab}u^b=0 ;\nonumber \\
&\dot{\mu}_A+\mu_B C^B_{Ab}u^b-\rho^i_A\pd{\cf}{x^i}+\pd{\cf}{u^c} C^c_{Ab}u^b=0.
\end{eqnarray}
These equations are also obtained
in~\cite{maruskin:2007p3969,MaruskinBloch:2010p2554} for the case
$E=TM$.


\begin{remark}
For simplicity we are considering only normal extremals. For abnormal extremals
we just have to consider the Hamiltonian function to be  $H=\mu_au^a$ and solve
the same equations, i.e.
\begin{equation}
\mu_a=0; \qquad
\dot{x}^i=\rho^i_au^ a;  \qquad
\mu_B C^B_{ab}u^b=0; \qquad
\dot{\mu}_A+\mu_B C^B_{Ab}u^b=0. 
\end{equation}
with $(\mu_A(t))\neq (0)$ for all $t$.
\end{remark}

\begin{remark}
In the particular case when $\D=E$ and $\sigma=\id_E$ we recover the
Euler-Lagrange equations on the Lie algebroid $E$ for the Lagrangian
$L=\cf$. Also when $\D\subsetneq E$ we get the so-called vakonomic
equations for the Lagrangian $L=\cf$ (see~\cite{Iglesias:2007p4370}) 
\end{remark}

\section{Dynamic Optimal Control}
\label{sec:dynam-optim-contr}
In the dynamic problem, we can control directly the motion on a nonholonomic
system, with the exception of the constraint forces, of course. For instance,
we can consider the equations of motion to be  $\delta L(z)|_\D=u$, with $u\in \D^*$
the control variables representing the external (generalized) forces acting on
the system. Another possibility would be to consider systems on which the
accelerations are the control variables. 

In both kinds of problems the state space is the manifold $\D$ and the control bundle 
$\map{\pi}{B}{\D}$ is 
\begin{equation}
  \label{eq:controlbundle}
B=\set{(z,\nu)\in\TDD\times\D^*}{z\in\mathrm{Adm}(E)\quad\mathrm{and}\quad
  \delta L(z)|_\D=\nu}.
\end{equation}
An element $(z,\nu)$ of $B$ is of the form $z=(a,a,v)$ with $a\in\D$ and
$v\in T_a\D$ and where $\nu$ is determined by the equations $\langle \nu,b\rangle=\langle\delta
L(z),b\rangle$ for every $b\in\D$.  

When we consider the forces as control variables, since we are assuming that the constrained
Lagrangian system is regular, we can identify $B$ with
$\map{\mathrm{pr}_1}{\D\oplus\D^*}{\D}$, via $(a,a,v;\nu)\equiv (a,\nu)$, because the vector
$v$ is determined 
by the point $a$ and the equation $\delta L(z)|_\D=u$.

When we consider the accelerations as controls we can identify $B$ with
$\TDD\cap\mathrm{Adm}(E)\to\D$, via $(a,a,v;\nu)\equiv (a,a,v)$ because $\nu$ is
determined by $\nu=\delta L(z)|_{\D}$. 

From a formal point of view both problems are equivalent, since the relation
between them is one-to-one and thus it is possible to use the optimal solution
written in terms of 
accelerations to determine the optimal forces and viceversa. In other words,
they are related by a feedback transformation.  

However, from the practical point of view the second problem produces simpler
expressions. Therefore, we can identify $B$ with $\TDD\cap\mathrm{Adm}(E)$
and take coordinates $(x^i,y^a,v^a)$ where $v^a$ are the acceleration
coordinates, i.e. our control variables. 

On this set we also need to
specify a control system where the optimization will be built. Such a system is
specified by giving a section $\sigma:B\to \TED$, i.e. the resulting system must
always define an admissible velocity and acceleration.  The Lie algebroid
relevant for this case is the $E$-tangent to $\D$. An element
of $\TED$ is of the form $z=(a,b,w)$ with $a\in\D$, $b\in E$, with $\tau(a)=\tau(b)$ and
$w\in T_a\D$ with $\rho(b)=T\tau(w)$. Taking a local basis $\{e_a\}$ for $\D$, and
completing a local basis $\{e_a,e_A\}$ for $E$, we can write $a=y^ae_a$,
$b=z^ae_a+z^Ae_A$, and $w=(\rho^i_az^a+\rho^i_Az^A)\pd{}{x^i}+w^a\pd{}{y^a}$. By
taking coordinates $(x^i)$ in the base, we have coordinates
$(x^i,y^a,z^a,z^A,w^a)$ on $\TED$. A local basis of sections of $\TED\to\D$ is
$\{\X_a,\X_A,\V_a\}$ and the element $z$ can be written
$z=z^a\X_a(x,y)+z^A\X_A(a,y)+w^a\V_a(a,y)$, and
$\rho^1(z)=w=(\rho^i_az^a+\rho^i_Az^A)\pd{}{x^i}+w^a\pd{}{y^a}$. The corresponding
coordinates on the dual bundle $(\TED)^*$ will be denoted
$(x^i,y^a,\mu_a,\mu_A,\pi_a)$.

If we choose to control the accelerations of the system $\{u^a\} $, the
map $\sigma:B\to \TED$ is given by the natural inclusion
$\sigma(z)=z$, 
\[
\xymatrix{& \TED \ar[d]_\tau \ar[r]^{\rho^1}& T\D\ar[dl]^{\tau_\D}\\
           B\ar[ur]^{\sigma}\ar[r]& \D &
}
\]
which in coordinates corresponds to
$\sigma(x^i,y^a,u^a)=(x^i,y^a,y^a,0,u^a)$.

Given a cost function $\cf:B\to\mathbb{R}$ we take the Pontryagin Hamiltonian $H\in
C^\infty((\TED)^*\times B)$, defined as $ H(\mu,z)=\langle\mu,
\sigma(z)\rangle-\cf(z) $ which in coordinates is  
\[
H(x^i,y^a, \mu_a, \mu_A, \pi_a,u^a)=\mu_ay^a+\pi_au^a-\cf(x^i,y^a,u^a),
\]
where the control functions are the accelerations $u^a$.

From $\pd{H}{u^a}=0$, we get
\begin{equation}
  \label{eq:piaopt}
\pi_a=\pd{\cf}{u^a}.
\end{equation}
These equations determine the optimal submanifold $S_H$ in $(\TED)^*\times_{\D} B$, which
in this case is a section of the projection $(\TED)^*\times_{\D} B\to B$, locally given  by
the equations 
$u^a=u^a(x^i, y^a, \pi_a)$. On $S_H$, the equations of motion are the
following. From $\dot{x}=\rho\pd{H}{\mu}$, since in this case the base variables are
$(x,y)$, we get the original control system 
\begin{equation}
  \label{eq:controlsys_d}
\dot{x}^i=\rho^i_ay^a+\rho^i_A 0=\rho^i_ay^a  \qquad
\dot{y}^a=u^a.
\end{equation}
Also, the equations of motion for $\pi_a$, written as
$\dot{\pi}_a=-\pd{H}{y^a}=\pd{\cf}{y^a}-\mu_a$,
because all the structure functions involved vanish (i.e. $\V_a$ commute with
all the others). Therefore we get
\[
\mu_a=\pd{\cf}{y^a}-\frac{d}{dt}\left(\pd{\cf}{u^a}\right)=\pd{\cf}{y^a}-\dot{\pi}_a,
\]
which is the combination that appear in~\cite{maruskin:2007p3969}, equation (11), for the
case $E=TM$, under the notation $\kappa_J$ and without a justification. 

The equation of motion for $\mu_a$ is
\begin{equation}
  \label{eq:mua}
-\dot{\mu}_a 
=\rho^i_a\pd{H}{x^i}+\mu_cC^c_{ab}y^b+\mu_CC^C_{ab}y^b
=-\rho^i_a\pd{\cf}{x^i}+\mu_cC^c_{ab}y^b+\mu_CC^C_{ab}y^b,
\end{equation}
and the equation of motion for $\mu_A$ is
\begin{equation}
  \label{eq:muA}
-\dot{\mu}_A
=\rho^i_A\pd{H}{x^i}+\mu_cC^c_{Ab}y^b+\mu_CC^C_{Ab}y^b 
=-\rho^i_A\pd{\cf}{x^i}+\mu_cC^c_{Ab}y^b+\mu_CC^C_{Ab}y^b.
\end{equation}
These two last equations correspond to equations (18) and (19) obtained
in~\cite{maruskin:2007p3969} for the case $E=TM$. 

\begin{remark}
As in the previous case, we have considered only normal extremals. For abnormal
extremals we just have to take the cost function $\cf=0$ and solve the same
equations for $\cf=0$, that is
\begin{equation}
\mu_a=0; \quad
\pi_a=0; \quad
\dot{x}^i=\rho^i_ay^a;\quad
\mu_B C^B_{ab}y^b=0; \quad
\dot{\mu}_A+\mu_B C^B_{Ab}y^b=0.
\end{equation}
with $(\mu_A(t))\neq (0)$ for all $t$. Interestingly, we get exactly the same
equations (plus $\pi_a=0$) as in the kinematic case.
\end{remark}

\section{Examples}
\label{sec:examples}
First we discuss from the point of view of the theory of Lie algebroids an example 
studied in~\cite{MaruskinBloch:2010p2554}. The result is naturally analogous to
the results obtained there, but within the new framework the treatement of
dynamical control systems becomes much more natural. After that, we
 study a few other examples of systems relevant for their
applications and whose solutions, obtained in more involved ways, can
be found in the literature.

\subsection{Dynamic optimal control of the vertical rolling disc}
For such a system, the state space manifold corresponds to  $M=\mathbb{R}^2\times
S^1\times S^1$, and we will use the coordinates $x=(x^1,x^2,x^3,x^4)$, where
$x^1=x,\;x^2=y,\;x^3= \theta,\;\mbox{and } x^4= \phi$.

The rolling without slipping condition of the motion on the plane leads to a
pair of nonholonomic constraints
$$
\dot x^1-\cos ( x^4 )\dot x^3=0\;\; \mathrm{\,\, and }\;\; \dot
  x^2-\sin ( x^4 )\dot x^3=0.
$$

We can define then a set of coordinates adapted to these constraints and write
a set of coordinates $\{ y \} $ for the new velocities. Thus the quasi-velocities
correspond to:
$$
y^1=\dot x^1-\cos( x^4)\dot x^3,\quad y^2=\dot x^2-\sin( x^4)\dot x^3,\quad
y^3=\dot  x^3,\quad y^4=\dot  x^4.
$$
Analogously the inverse transformation allows us to write:
$$
\dot x^1=y^1+\cos( x^4)y^3,\quad \dot x^2=y^2+\sin( x^4)y^3,\quad \dot
x^3=y^3,\quad \dot x^4=y^4.
$$

 The  local basis of sections of $TM$ determined by the quasi-velocities turns
 out to be:
$$
     e_1=\frac{\partial}{\partial x^1}, \quad e_2=\frac{\partial}{\partial x^2}, \quad
     e_3=\cos( x^4)\frac{\partial}{\partial x^1}+\sin( x^4)\frac{\partial}{\partial
       x^2}+\frac{\partial}{\partial x^3}, \quad
    e_4=\frac{\partial}{\partial x^4}.
$$
The Lie algebroid structure is the usual one for the tangent bundle. But in
the basis above, the anchor mapping is written as:
$$
    \rho(e_1)=\frac{\partial}{\partial x^1}, \quad
    \rho(e_2)= \frac{\partial}{\partial x^2}, \quad
    \rho(e_3)=\cos( x^4)\frac{\partial}{\partial x^1}+\sin( x^4)\frac{\partial}{\partial
      x^2}+\frac{\partial}{\partial x^3},\quad
    \rho(e_4)=\frac{\partial}{\partial x^4}.
$$
The Lie algebra structure is obtained also as
$
[e_3,e_4]=\sin (x^4)e_1-\cos (x^4)e_2,
$
all other elements being zero. We can read then the Hamel symbols 
$\gamma_{\alpha\beta}^\epsilon$  in \cite{Carinena:2007p2820,maruskin:2007p3969} .


In what regards the control part, we are considering a situation where we
control the external forces in the directions of the admissible velocities.
Thus, as the velocities on the constrained system are of the form
$$
y^1=y^2=0,\;\; y^3=\dot x^3\qquand y^4=\dot x^4,
$$
the natural coordinates are  $\mathrm{a}=(x^1,x^2,
x^3,x^4,y^3,y^4)=(x,y^a)$. The control bundle $B$ becomes thus
$\D\oplus \D^*$ and we 
 take as coordinates $(x^i, y^3, y^4,u_3, u_4)$, where $u_3= \frac
 32\dot y^3$ and $u_4= \frac 14\dot y^4$ for $u_a=(\delta L)_a$. 

The cost function corresponds to
$
\kappa(\mathrm{a},u)=\frac 12(u_3^2+u_4^2)
$
and the control system is defined as:
$$\dot x^1=\cos( x^4)\dot x^3, \quad \dot x^2=\sin( x^4)\dot x^3,\quad
u_3=\frac 32\ddot{ x}^3,\quad u_4= \frac 14\ddot{ x}^4.
$$
The Pontryagin Hamiltonian  $\;H\in C^\infty((\TED)^*\times_{\D} B)$ corresponds
now to
\[
H(\mathrm{a},p,u)=\langle p, \sigma_\mathrm{a}(u)\rangle-\kappa(\mathrm{a},u)
=\mu_Iy^I+\pi_I\frac{u_I}{c_I}-\frac 12(u_3^2+u_4^2)
,
\]
with $c_3=3/2$ and $c_4=1/4$.

The Maximum principle is encoded as
\begin{center}
\begin{tabular}{lllll}
  $\displaystyle{\frac{\partial H}{\partial u_3}}=0$ & $\Leftrightarrow u_3=\frac 23\pi_3$ $\qquand$& &
  $\displaystyle{\frac{\partial H}{\partial u_4}}=0$ & $\Leftrightarrow u_4=4\pi_4$,
\end{tabular}
\end{center}
and the Pontryagin equations (optimal dynamical control equations) correspond to:
\begin{equation*}
 \dot x^1=\cos( x^4)y^3, \qquad  \dot x^2=\sin( x^4)y^3, \quad
\dot x^3= y^3, \qquad \dot  x^4=y^4
\end{equation*}
\begin{equation*}
\dot y^3=\frac 23 u_3, \qquad \dot y^4=4u_4, \quad
 \dot\pi_3=-\mu_3, \qquad \dot\pi_4=-\mu_4,
  \end{equation*}
  \begin{equation*}
\dot \mu_1=0, \qquad \dot\mu_2=0,
 \end{equation*}
\begin{equation*}
\dot\mu_3=[-\mu_1\sin( x^4)+\mu_2\cos( x^4)]y^4, \qquad
\dot\mu_4=[\mu_1\sin( x^4)-\mu_2\cos( x^4)]y^3.
\end{equation*}
Since $y^3=\dot x^3$, $y^4=\dot x^4$ and $\dot y^I=u_I$, then
$\mu^3=-\frac{9}{4}{\dddot{x}}^3\,\mbox{and } \mu^4=-\frac{1}{16}\dddot x^4$.
Thus, we can reduce the set of equations to:
\begin{align*}
  \dot x^1&=\cos( x^4)\dot x^3,\quad                                         
\dot x^2=\sin( x^4)\dot x^4,\\
  \ddddot{x}^3&=\frac{4}{9}[\mu_1\sin( x^4)-\mu_2\cos( x^4)]\dot x^4,  \\
 {\ddddot x}^4&=16[-\mu_1\sin( x^4)+\mu_2\cos( x^4)]\dot x^3,
\end{align*}
where $\mu_1,\mu_2$ are constants.

\subsection{Optimal control problems of rotational motion of the  free rigid body} 

Consider the problem of rotational motion of the free rigid body. As
configuration manifold  we take the Lie group $SO(3)$ and choose the type-I
Euler angles $(x^1,x^2, x^3)$ as
  local coordinate system. We consider the canonical  Lie algebroid
 structure of  the tangent bundle $TSO(3)$, whose anchor map is
 $\rho=id_{TSO(3)}$. Let $\{e_1,e_2,e_3\}$ be the set of sections
 for the bundle
\begin{align*}
e_1&= \mathrm{sec}(x^2) \sin (x^3) \frac{\partial}{\partial x^1}+ \cos
(x^3) \frac{\partial}{\partial x^2}+ \tan (x^2) \sin (x^3)\frac{\partial}{\partial x^3},
\nonumber \\
e_2&= \mathrm{sec}(x^2) \cos (x^3)
\frac{\partial}{\partial x^1}- \sin (x^3) \frac{\partial}{\partial
x^2}+\tan (x^2) \cos (x^3) \frac{\partial}{\partial x^3},
\\
e_3&=
\frac{\partial}{\partial x^3},  
\end{align*}
whose Lie algebra structure is determined by the relations  $[e_1,e_2]=e_3$,
$[e_2,e_3]=e_1$, $[e_3,e_1]=e_2$. The anchor and the Lie bracket are locally
determined by the functions
$$
\begin{array}{lll}
\rho^1_1=\mathrm{sec}(x^2) \sin (x^3), & \rho^2_1=\cos (x^3), &   \rho^3_1=\tan
(x^2) \sin (x^3), \\
\rho^1_2=\mathrm{sec}(x^2) \cos (x^3), & \rho^2_2=- \sin (x^3),  &
\rho^3_2=\tan (x^2) \cos (x^3), \\
\rho^1_3=0,   & \rho^2_3=0, & \rho^3_3=1;
\end{array}
$$
and
$C_{12}^3=-C_{21}^3=C_{23}^1=-C_{32}^1=C_{31}^2=-C_{13}^2=1$.

We consider now the rigid body  subject to the constraint
$
\dot x^1\cos (x^2)\sin (x^3)+\dot x^2\cos (x^3)=0
$.
This implies that the constraint distribution $\D$ is
 the $2$-dimensional subbundle of $TSO(3)$ generated by $e_2$ and $e_3$.

\subsubsection{Constrained kinematic problem}

Let us study the kinematic constrained system, i.e., the system with admissible
velocities belonging to the subbundle $\D\subset TSO(3)$ defined by the condition
$
y^1=0
$.
Thus the system is
 \begin{equation}
   \label{eq:xs3vv}
\dot{x}^1=\mathrm{sec}(x^2)\cos (x^3) u^2, \quad
\dot{x}^2=-\sin (x^3) u^2, \quad
\dot{x}^3=\tan (x^2)\cos (x^3) u^2+u^3,
 \end{equation}

The cost function corresponds to the energy provided by
the controls
$
k(x^i,u^a)=\frac{1}{2}\left[I_{2}(u^2)^2+I_{3}(u^3)^2\right]
$.
The Hamiltonian in this case is written as
$$
H=\mu_2u^2+\mu_3u^3-\frac{1}{2}\left[I_{2}(u^2)^2+I_{3}(u^3)^2\right].$$
Optimality conditions defining the submanifold $S_{H}$ given in (\ref{eq:optimal}) are
$
\mu_2=I_2u^2;\mu_3=I_3u^3
$.
Using the representation $u^2=u^2(x, \mu)$ and $u^3=u^3(x, \mu)$ for $S_{H}$, the
equations (\ref{eq:dotm-=-leftrh}) become then the control system
(\ref{eq:xs3vv}) together with
$$
\dot{\mu}_2 +\mu_1u^3=0, \quad
\dot{\mu}_3 -\mu_1u^2=0, \quad
\dot{\mu}_1+(I_3-I_2)u^2u^3=0.
$$

In the case of the completely symmetric rigid body  we get
$
\dot{\mu}_2 +\mu_1u^3=0; \quad
\dot{\mu}_3 -\mu_1u^2=0$ and 
 $\dot{\mu}_1=0
$,
which are equivalent to the equations obtained by  Sastry and Montgomery in
 \cite{Sastry:1993p4035}.


\subsubsection{Constrained dynamic problem}
Finally let us study the case of dynamic control for the constrained
system. From the geometrical point of view, the control bundle $B$
corresponds to $\TDD\cap \rm{Adm}(TSO(3))$ and the system can be
described as
 \begin{align}
   \label{eq:xs3v}
&\dot{x}^1=\mathrm{sec}(x^2)\cos (x^3) y^2, \quad
\dot{x}^2=-\sin (x^3) y^2, \\
&\dot{x}^3=\tan (x^2)\cos (x^3) y^2+y^3, \quad
\dot{y}^2=u^2,\;\dot{y}^3=u^3
 \end{align}

We consider now as cost function, the restriction to $\D$ of the
cost function
$$
\begin{array}{l}k(x^i,y^a,u^a)=\frac{1}{2}\left[(I_2)^2(u^2)^2+(I_3)^2(u^3)^2
+(I_3-I_2)^2(y^2)^2(y^3)^2\right]\end{array},$$
where we assume that the control functions are the components of the admissible angular
accelerations of our system  $u^2=v^2$ and $u^3=v^3$.

Optimality condition leads to the submanifold $W$ defined as
$
\pi_2=(I_2)^2u^2 \quad
\pi_3=(I_3)^2u^3
$. Then $W$ is defined by specifying $u^2=u^2(x, y, \pi)$ and $u^3=u^3(x, y,
\pi)$. The motion on $W$ corresponds then to the control system (\ref{eq:xs3v}) and
 \begin{align*}
&\dot{\pi}_2=\frac{(M_1)^2}{y^2}-\mu_2,  \quad
\dot{\pi}_3=\frac{(M_1)^2}{y^3}-\mu_3,  \quad
\dot{\mu}_2 +\mu_1y^3=0, \quad
\dot{\mu}_3 -\mu_1y^2=0, \\
&\dot{\mu}_1+\mu_3y^2-\mu_2y^3=0,
 \end{align*}
where $M_1=(I_3-I_2)y^2y^3$ is a torque on $\D$.

In the case of the completely symmetric rigid body
we obtain
 $$
\dot{\pi}_2=-\mu_2,  \quad
\dot{\pi}_3=-\mu_3,  \quad
\dot{\mu}_2 +\mu_1y^3=0, \quad
\dot{\mu}_3 -\mu_1y^2=0, \quad
\dot{\mu}_1+\mu_3y^2-\mu_2y^3=0.
$$
This system gives the following equations obtained by Crouch and Silva Leite in
\cite{Crouch:1995p3341}, Ex. 6.4, Case II
$$
\dddot{y}^2 -\mu_1y^3=0, \quad
\dddot{y}^3+\mu_1y^2=0, \quad
\dot{\mu}_1-\ddot{y}^3y^2+\ddot{y}^2y^3=0.
$$

\subsection{Systems with symmetry and constraints: quasi-coordinates for the
  Atiyah algebroid}

Consider a ball rolling without sliding on a fixed table (see Example 8.12 in
\cite{Cortes:2009p2371}). The configuration space is $Q=\mathbb{R}^2\times SO(3)$, where
$SO(3)$ is parameterized by the Eulerian angles $\theta,\phi$ and $\psi$. In
quasi-coordinates $(x,y,\theta,\phi,\psi,\dot x, \dot y, \omega_x,\omega_y,\omega_z)$ the energy may be
expressed by
    $T=\frac{1}{2}[\dot x^2+\dot y^2 + k^2(\omega_x^2+\omega_y^2+\omega_z^2)]$,
where $\omega_x,\omega_y$ and $\omega_z$ are the components of the angular
velocity of the ball.

The system is invariant under $SO(3)$ transformations, and thus it is natural
to consider the corresponding formulation on the Atiyah algebroid $E=TQ/SO(3)\equiv
T\mathbb{R}^2\times\mathbb{R}^3$.  On that system we must still implement the
nonholonomic constraint arising from the rolling-without-sliding conditions
$\dot x^1-r\omega_2=0$ and $\dot x^2+r\omega_1=0$.

For the configuration space we can choose coordinates
$
M=Q/G=\mathbb{R}^2 \ni \mathrm{x}=(x^1,x^2)$, with $x^1=x$ and $x^2=y$.
In what regards the fiber, we can choose thus a transformation mapping the set
of fiber coordinates $\{\dot x^1, \dot x^2, \omega_3,\omega_1,\omega_2\}$ onto a
new set $\{ y^\alpha\}$. These quasi-velocities become then
$ y^i=\dot x^i,\quad y^3=\omega_3,\quad y^4=\dot x^1-r\omega_2,\quad y^5=\dot
x^2+r\omega_1$.

Analogously we can consider the inverse transformation.
Thus the original velocities can be written in terms of the quasi-velocities as
$$
\dot x^i=y^i,\quad \omega_3=y^3,\quad
\omega_1=-\frac{1}{r}y^2+\frac{1}{r}y^5,\quad\omega_2=\frac{1}{r}y^1-\frac{1}{r}y^4
$$


The local basis of sections of $E$ determined by the quasi-velocities turns out
to be
$$
f_1=e'_1+\frac{1}{r}e'_4,\quad
  f_2=e'_2-\frac{1}{r}e'_3,\quad
f_3=e'_5,\quad
f_4=-\frac{1}{r}e'_4,\quad
f_5=\frac{1}{r}e'_3
$$
where $\{ e'_1, e'_2, \cdots , e'_5\}$ is the local basis defined in \cite{Cortes:2009p2371}, page 36.

With respect to this basis, the structure constants and the anchor mapping  of
the Lie algebroid structure become
$$
  [f_2,f_1]=[f_1,f_5]=[f_4,f_2]=[f_5,f_4]=\frac{1}{r^2}f_{3}, 
$$
$$[f_3,f_1]=[f_4,f_3]=f_5 [f_2,f_3]=[f_3,f_5]=f_4,
$$
 $$
\rho(f_1)=\partial_x,\quad \rho(f_2)=\partial_y, 
$$
the remaining elements being zero.

The set of admissible velocities becomes thus the fiber of the distribution
$\D$, which corresponds to
$
\D=\{(x^i,y^\alpha)\in E\mid y^4=y^5=0\}
$.
The coordinates for these points are therefore $\mathrm{a}=(x^1,x^2,y^1,y^2,y^3)=
(x^i,y^\alpha)$, where in terms of the original set of coordinates these correspond
to $y^1=\dot x^1=u^1$,  $y^2=\dot x^2=u^2$, $y^3=\omega_3=u^3$ and $y^4=0=y^5$.

The dynamical system on the algebroid  is defined by a Lagrangian function on
$E$, which can be written in terms of the velocities as $L(x,y,\dot x, \dot y, \omega_1,\omega_2,\omega_3)
    =\frac{1}{2}[\dot x^2+\dot y^2 +    k^2(\omega_1^2+\omega_2^2+\omega_3^2)]$, 
and in terms of the quasi-velocities
  \begin{align}
    \label{eq:L2}
  L(x^i,y^\alpha)
    &=\frac{1}{2}\left [ (y^1)^2+(y^2)^2 +\right .\nonumber  \\
  &\left . k^2r^{-2}((y^1)^2+(y^2)^2+(y^4)^2+(y^5)^2-y^2y^5-y^1y^4)+
    k^2(y^3)^2\right ].
  \end{align}

\subsubsection{Kinematic Control Problem}
Consider the following problem: determine the minimal value among the set
of admissible solutions $a:\mathbb{R}\to E$, of the controlled  Euler-Lagrange
equations of the form $\delta L(a(t))=0, \quad a(t)\in \D$,
where the cost fuction is
$\kappa(x^i,u^a)=
    \frac{1}{2}\left\{(u^1)^2+(u^2)^2 +
      \frac{k^2}{r^2}\left[(u^2)^2+(u^1)^2\right]+k^2(u^3)^2\right\}$.
The control bundle is  $ \D$, and the section we consider $\sigma: (x^i, u^a)\in \D\to
 (x^i, u^a,0,u^a)\in E$ is the canonical inclusion. Please notice that we use the notation $u^a$ to denote the
 elements of the fiber of $\D$ when considered as the control bundle. 
%

The  Pontryagin Hamiltonian is then written as a function
$H\in C^\infty(E^*\times_{\mathbb{R}^2} \D)$.  The optimality conditions of the Maximum
principle on this function imply

 \begin{equation}
   \label{eq:extremal}
   \frac{\partial H}{\partial u^a}=0 \Leftrightarrow
   \mu_a=c_au^a,
 \end{equation}
 where $c_1=c_2=1+k^2/r^2$ and $c_3=k^2$.

If we write the set of Pontryagin equations we see:
\begin{align*}
\dot x^i&=&y^i, \qquad \qquad \qquad\,\,\, 
\dot \mu_1&=&\mu_3\frac{\mu_2}{c_2r^2}+\mu_5\frac {\mu_3}{c_3},\quad
\dot \mu_2&=& -\mu_4\frac{\mu_3}{c_3}-\mu_3\frac{\mu_1}{c_1r^2},\\
\dot \mu_3&=& -\mu_5\frac{\mu_1}{c_1}+\mu_4\frac{\mu_2}{c_2}, \quad
\dot \mu_4&=&-\mu_3\frac{\mu_2}{c_2r^2}-\mu_5\frac{\mu_3}{c_3},\quad
\dot \mu_5&=&\mu_4\frac{\mu_3}{c_3}+\mu_3\frac{\mu_1}{c_1r^2}.
\end{align*}
Thus we can use (\ref{eq:extremal}) and the above equations to define the
resulting system on $\D$:
\begin{equation*}
\ddot x^1=\frac{1}{c_1}(d_2\omega_3-\dot x^2\omega_3),\quad
\ddot x^2=\frac{1}{c_2}(\dot x^1\omega_3-d_1\omega_3),\quad
\dot \omega_3=\frac{1}{c_3}(d_1\dot x^2 -d_2\dot x^1),
\end{equation*}
with $d_1,d_2$ constants and $c_1=c_2=1+k^2/r^2$ and $c_3=k^2$.\\

 \subsubsection{Dynamic Optimal Control Problem}
Let us consider now a different control problem, where we are able to
control the forces acting on the system, i.e. we consider a system
corresponding  to
$
(\delta L)_a=u_a
$,
where $L$ is defined as (\ref{eq:L2}) and $\delta$ represents the variational
derivative. In this case, for the Lagrangian given above, this implies that
$
\dot y^a=\displaystyle{\frac{u_a}{c_a}}, \quad a=1, 2, 3;
$
 where again $c_1=c_2=1+k^2/r^2$ and $c_3=k^2$ (see equations 1.9.13 in \cite{Bloch:2003p2763}).

The control system is thus defined as a section
$\;\sigma:(x^i, y^a, u_a)\in \D\oplus \D^*\mapsto (x^i, y^a, y^a,0,u_a/c_a)\in \TED$.
The cost function now is the energy provided by the control functions:
$
 \cf(x^i, u_a)=\frac{1}{2}\sum_{a}u_a^2
$.
As a result, the Pontryagin Hamiltonian $H\in C^\infty((\TED)^*\times_{\D} B)$
reads now
$
H(x^i, y^a,\mu_\alpha,\pi_a, u_a)=\mu_a y^a+\pi_a \frac{u_a}{c_a}-\frac{1}{2}\sum_{a}u_a^2
$.
The maximum principle applied to this function results
$$
\frac{\partial H}{\partial u_a}=0 \Leftrightarrow u_a=\frac{\pi_a}{c_a}, \,\,{\rm with}\,\, a=1,2,3.
$$
Then the optimal manifold correponds to this submanifold of
$(\TED)^*\times_{\D} B$.

The  Pontryagin equations on $(\TED)^*\times_{\D} B$ are:

$$
\dot x^i= y^i, \qquad
\dot y^a=\frac{u_a}{c_a}, \qquad 
\dot \pi_a=-\mu_a 
$$
$$
\dot \mu_1=\mu_3\frac{y^2}{r^2}+\mu_5y^3,
\qquad
\dot \mu_2= -\mu_4y^3-\mu_3\frac{y^1}{r^2},\qquad 
\dot \mu_3=-\mu_5y^1+\mu_4y^2, 
$$
$$
\dot \mu_4=-\dot\mu_1, \qquad
\dot \mu_5=-\dot\mu_2.
$$
But if we restrict them to the optimized submanifold we obtain the reduced system:
\begin{align*}
  &\ddddot x^1=[\frac{c_3}{c_1r}]^2\dot x^2\ddot\omega_3-\omega_3\dddot
  x^2-\frac{e_2}{c_1^2}\omega_3,\\
  &\ddddot x^2=-[\frac{c_3}{c_2r}]^2\dot x^1\ddot\omega_3+\omega_3\dddot
  x^1+\frac{e_1}{c_2^2}\omega_3,\\
  &\dddot {\omega_3}=[\frac{c_2}{c_3}]^2\dot x^1\dddot x^2-[\frac{c_1}{c_3}]^2\dot x^2\dddot x^1
  +\frac{e_2}{c_3^2}\dot x^1-\frac{e_1}{c_3^2}\dot x^2,
\end{align*}
where $c_1=c_2=1+k^2/r^2,\;c_3=k^2$ and $e_1,e_2$ are arbitrary constants.

\section*{Acknowledgments}

Partial financial support from MEC-DGI (Spain) grants
 MTM2009-11154, MTM2009-08166-E, Fis2009-12648, the portuguese grant
 PTDC/MAT/099880/2008 and Aragon's E24/1 (DGA)  is acknowledged.

\end{document}